\documentclass{article}

\usepackage{amsmath}
\usepackage{enumitem}
\usepackage{amsfonts}
\usepackage{amssymb}
\usepackage{graphicx}
\usepackage{mathrsfs}
\usepackage{esint}
\usepackage{upref,amsthm,amsxtra,exscale}
\usepackage{cite}
\usepackage[colorlinks=true,urlcolor=blue,
citecolor=red,linkcolor=blue,linktocpage,pdfpagelabels,
bookmarksnumbered,bookmarksopen]{hyperref}
\usepackage{cleveref}
\usepackage[cm]{fullpage}
\usepackage{subcaption}
\usepackage{caption}

\newtheorem{theorem}{Theorem}[section]
\newtheorem{corollary}[theorem]{Corollary}
\newtheorem{remark}[theorem]{Remark}
\newtheorem{lemma}[theorem]{Lemma}
\newtheorem{proposition}[theorem]{Proposition}

\newtheorem{example}[theorem]{Example}
\numberwithin{equation}{section}

\title{Multiplicity of sign-changing solutions to higher-order $Q$-curvature equations}
\author{Juan Carlos Fernández\thanks{Supported by UNAM-DGAPA-PAPIIT IN102027 grant.}}
\date{ }

\begin{document}
\maketitle

\begin{abstract}
In this paper, we employ the dual cone method to establish the existence of infinitely many sign-changing solutions to the $Q$-curvature problem involving $GJMS$-operators of arbitrary order $m$ on closed Riemannian manifolds of dimension $n > 2m$, under the assumptions that the operators are coercive and positivity-preserving and that the manifold admits appropriate symmetries. To demonstrate the applicability of our framework, we provide explicit examples on both closed Einstein manifolds and closed special Einstein products. Notably, one of the major contributions of this work is deriving explicit conditions that guarantee the coercivity of $GJMS$-operators on closed special Einstein products.

\medskip


\medskip


\end{abstract}

\section{Introduction}

Partial differential equations arising from geometric problems have attracted considerable attention over the past
several decades.  Among the central problems in this area is the prescribed $Q$-curvature problem, which leads to critical
Sobolev exponent equations involving conformally covariant differential operators. More precisely, let $(M,g)$ be a closed (compact and without boundary) Riemannian manifold of dimension $n$, and let $m\in\mathbb{N}$ be such that $2m<n$. In 1992,
Graham, Jenne, Mason and Sparling introduced in \cite{GrahamJenneMasonSparling1992}  a linear partial differential operator $\mathscr{P}_g^m:C^\infty(M)\to C^\infty(M)$ with leading term $(-\Delta_g)^m$, where $\Delta_g = \text{div}_g(\nabla\cdot)$ denotes the Laplace-Beltrami operator on $(M,g)$, satisfying the following properties:
\begin{enumerate}[label=($P_1$),ref=$(P_1)$]
	\item \label{Propoerty:Pg naturality} $P_g$ is natural, i.e., for every diffeomorphism $\varphi:M\rightarrow M$,  $\varphi^\ast P_{g}=P_{\varphi^\ast g}\circ\varphi^\ast$, where ``$\!\!\!\phantom{K}^\ast\;$'' denotes the pullback of a tensor. 
    \end{enumerate}
  \begin{enumerate}[label=($P_2$),ref=$(P_2)$]
	\item \label{Property:Pg SelfAdjoint} $P_{g}$ is self-adjoint with respect to the $L^2$-scalar product.
    \end{enumerate}
	\begin{enumerate}[label=($P_3$),ref=$(P_3)$]
    \item\label{Property:Pg ConformalInvaiance} $P_g$ is conformally invariant, that is, given any function $\omega\in\mathcal{C}^\infty(M)$, if we define the conformal metric $\widetilde{g}=e^{2\omega} g$, then the following identity holds true:
	\begin{equation*}\label{Eq:GJMS-conformal invariance}
		P_{\widetilde{g}}(f)=e^{-\frac{N+2m}{2}\omega}P_{g}\left(e^{\frac{N-2m}{2}\omega}f\right), \text{ for all }f\in\mathcal{C}^\infty(M).   
	\end{equation*}
\end{enumerate}
The $Q$-curvature associated to $\mathscr{P}_g^m$ is defined by the function
\[
Q_g:= \frac{2}{n-2m}\mathscr{P}_g^m(1)\in C^\infty(M).
\]
The problem of finding a conformal metric $\widetilde{g}$ on $M$ for which $Q_{\tilde{g}}$ is equal to a constant $\mu\in\mathbb{R}$ is known as the prescribed constant $Q$-curvature problem. For $u\in C^\infty(M)$ is positive, if we write $\widetilde{g}=u^{\frac{4}{n-2m}}g$, then, by property \ref{Property:Pg ConformalInvaiance}, this problem is equivalent to solving the PDE
\begin{equation}\label{MainEquation:Q-curvatureEquation}
\mathscr{P}_g^m u = \mu \vert u\vert^{2_{m,n}^\ast-2}u\qquad \text{on }M,
\end{equation}
where $2_{m,n}^\ast:=\frac{2n}{n-2m}$ denotes the critical Sobolev exponent of order $m$. Equation \eqref{MainEquation:Q-curvatureEquation} is known as the prescribed constant $Q$-curvature equation and  has attracted substantial interest in both differential geometry and the
theory of partial differential equations. 
The case $m=1$ corresponds to the famous Yamabe problem, where $\mathscr{P}_g^1=-\Delta_g  \frac{n-2}{4(n-1)}R_g$ is known as the conformal Laplacian, $R_g$ denotes the scalar curvature and, up to a dimensional constant, it is equal to $Q_g^1$, the $Q$-curvature in this case. For $m=2$, the operator $\mathscr{P}_g^2$ coincides with the Paneitz-Branson operator and $Q_g^2$ is the usual $Q$-curvature \cite{Paneitz1983,Branson1985}. 

The existence and multiplicity theory of positive solutions to equation \eqref{MainEquation:Q-curvatureEquation} has a long history. Concerning the Yamabe problem, corresponding to $m=1$, we refer the reader to the book \cite{Aubin1998}, the survey \cite{BrendleSurvey2008} and the references therein. In case of the Paneitz-Branson equation, the existence and multiplicity of solutions has been studied by several authors; see, for instance \cite{DjadliHebeyLedoux2000,QingRaske2006, Saintier2006,GurskyMalchiodi2015,HangYang2015,HangYang2016,GurskyHangLin2016,BettiolPiccioneSire2021,JulioBatallaPetean2024},  although this list is far from exhaustive. Considerable progress has also been made on the compactness of positive solutions to the Paneitz-Branson equation; see, for example, \cite{HebeyRobert2004,Malchiodi2006,Lin1998,WeiZhao2013,LiXiong2019}. Only more recently has the general case $m\geq 3$ been successfully addressed. Concerning the existence and multiplicity of solutions to \eqref{MainEquation:Q-curvatureEquation}, see, for instance, \cite{Robert2011,CaseMalchiodi2024,MazumdarVetois2024,AndradeCasePiccioneWei2026,MazumdarNdiaye2026}; meanwhile,  the compactness of the set of positive solutions has been studied in \cite{MazumdarPremoselli2025,GongKimWei2026}.

The situation for sign-changing solutions to \eqref{MainEquation:Q-curvatureEquation} is rather different. The existence and multiplicity of sign-changing
solutions to the Yamabe problem have been studied extensively; see, for instance, \cite{Ding1986,AmmannHumbert2006,ElSayed2014,RobertVetois2015,ClappFernandez2017,DelPinoMussoPacardPistoia2011,CorroFernandezPerales2023,FernandezPalmasTorresOrozco2024,ClappSaldanaSzulkin2021,ClappFayaSaldana2024,Stefaniak2026}, and, again, this list is not exhaustive. By contrast, substantially fewer results are available for the general problem  $m\geq 2$, mainly due to the lack of the maximum principle fo higher-order operators. For $\mu<0$ in equation \eqref{MainEquation:Q-curvatureEquation}, some progress has been achieved by adapting the strategy based on the second Yamabe
invariant developed in \cite{AmmannHumbert2006}. For the Paneitz-Branson operator, see \cite{BenaliliBoughazi2016,AzaizBoughazi2020,BoughaziTouati2024}, while the general case $m\geq 2$ is considered in \cite{AzaizBoughazi2025}. When $\mu>0$, the situation is more delicate, since the variational quantities associated with higher eigenvalues of the
nonlinear problem \eqref{MainEquation:Q-curvatureEquation} are substantially more difficult to attain, and in some settings attainment may fail.  One approach
to overcoming the lack of a maximum principle is to restrict attention to Einstein manifolds, for which the operator $\mathscr{P}_g^m$ admits a particularly useful factorization as a product of second-order conformally covariant operators  \cite{Gover2006}. In this direction, Julio-Batalla and Petean proved the existence nfinitely many sign-changing solutions to the Paneitz-Branson equation when $(M,g)$ admits a proper isoparametric function. As a particular case, they obtained such
a multiplicity result for Einstein manifolds with positive scalar curvature admitting a cohomogeneity-one action. In \cite{FernandezPalmasTorresOrozco2024}, an analogous multiplicity result was established for the same class of manifolds under the existence of group actions with positive-dimensional orbits. The first result relies on compactness properties of positive solutions with isoparametric
symmetries, whereas the second relies on the nonexistence of nonconstant positive solutions to the Paneitz-Branson
equation \cite{Vetois2024}. 

In the general case $m\geq 2$, very little is known about the existence of sign-changing solutions to the prescribed $Q$-curvature equation. When $\mu<0$, Azaiz and Boughazi in \cite{AzaizBoughazi2025}
 roved the existence of sign-changing solutions under the assumptions of positive $Q$-curvature and negativity of the first eigenvalue of $\mathscr{P}_g^m$. For $\mu>0$, as far as we know, the only available result concerning the existence and multiplicity of sign-changing solutions is due to Bartsch, Schneider and Weth \cite{BartschSchneiderWeth2004},  in the particular case in which the underlying manifold is the round sphere; their multiplicity result is obtained by imposing group actions with
positive-dimensional orbits. For existence results concerning nodal solutions to $Q$-curvature-type problems with a small perturbation of the nonlinearity, see also  \cite{OunaneTahri2026}.

\bigskip

Motivated by the last results, in this paper, we study the existence and multiplicity of sign-changing solutions to the general problem \eqref{MainEquation:Q-curvatureEquation} with $m\geq 2$ and $\mu>0$. To do this, we apply the dual cone method, presented in an abstract framework in \cite{Weth2006}, which has been proven useful for establishing  the existence of sign-changing solutions to equations involving critical and subcritical nonlinearities; see, for instance, \cite{BartschLiuWeth2004,ClappWeth2005,Weth2006,ClappPacella2008,ClappFernandez2017,CorroFernandezPerales2023}. We apply this method to closed Riemannian manifolds for which the operator $\mathscr{P}_g^m$ is coercive and preserves positivity, using symmetries as in \cite{FernandezPalmasTorresOrozco2024} to handle the lack of compactness due to the critical Sobolev exponent $2_{m,n}^\ast$. In order to state our main results, we establish these concepts and related concepts. 

\medskip

We say that the operator $\mathscr{P}_g^m$ is coercive if there exists $C>0$ such that
\[
\int_M u\mathscr{P}_g u\; dV_g \geq C \Vert u\Vert_{L_g^2(M)}^2,\qquad\forall\; u\in C^\infty(M).
\]
On the other hand, we say that $\mathscr{P}_g^m$ satisfies the positivity preserving property if
\begin{equation}\label{Property:PPP}
\tag{PPP}\text{for any } u\in C^\infty(M) \text{ such that }\mathscr{P}_g^m u\geq 0,\text{ it follows that }u\geq0.
\end{equation}
If $u$ is a nontrivial solution to \eqref{MainEquation:Q-curvatureEquation}, its energy is defined as the quantity
\[
\mu \frac{m}{n}\int_M \vert u\vert^{2_{m,n}^\ast}\; dV_g
\]

Now, we briefly recall the terminology regarding group actions on manifolds. If $\text{Isom}(M,g)$ denotes the group of isometries of $(M,g)$, there is an action by isometries on $M$  defined by $(\gamma, x)\mapsto \gamma x:=\gamma(x), x\in M,\gamma\in \text{Isom}(M,g)$. Given a closed subgroup $\Gamma\subset \text{Isom}(M,g)$, we say that a function $u:M\to\mathbb{R}$ is \emph{$\Gamma$-invariant} if  $u(\gamma x)=u(x)$ for every $x\in M$ and every $\gamma\in \Gamma$.  Given $x\in M$, define the $\Gamma$-orbit of $x$ as the set
\[
\Gamma x:= \{ \gamma x\in M \;:\; \gamma\in\Gamma \}.
\]
Since $\Gamma$ is closed, the $\Gamma$-orbits are submanifolds of $M$ (see, for instance, \cite{AlexandrinoBettiol2015}). We can state our main result.

\begin{theorem}\label{Theorem:Main}
Let $\mu>0$, $m\in\mathbb{N}$, $(M,g)$ be a closed Riemannian manifold of dimension $n> 2m$ and $\Gamma$ be a closed subgroup of $\text{Isom}(M,g)$. If the $GJMS$-operator $\mathscr{P}_g^m$ is coercive, satisfies  \eqref{Property:PPP}, and 
\begin{equation}\label{OrbitsPositiveDimensional}
1\leq \dim\Gamma x\leq n-1,\qquad\forall\; x\in M.
\end{equation}
then, the problem \eqref{MainEquation:Q-curvatureEquation} admits an infinite number of $\Gamma$-invariant solutions with increasing energies. One of them is positive, and the rest change sign.
\end{theorem}

We remark that condition \eqref{OrbitsPositiveDimensional} has two main purposes: the inequality $\dim\Gamma x<m$ guaranties that the action is not transitive, implying that underlying symmetric Sobolev spaces are infinite-dimensional; meanwhile, the inequality $1\leq \dim\Gamma x$ is needed to restore compactness of the Sobolev embeddings at the critical level. See  Section \ref{Section:VariationalSettingSymmetries} for the details.

\bigskip

The principal examples of closed Riemannian manifolds satisfying both coercivity and property \eqref{Property:PPP} are Einstein manifolds with positive scalar curvature; see \cite{Robert2011}. A fundamental example that also satisfies the symmetry condition \eqref{OrbitsPositiveDimensional} is the round sphere $\mathbb{S}^n$, endowed with the action considered in \cite{BartschSchneiderWeth2004}, namely $\Gamma=O(n_1)\times O(n_2)$ with $2\leq n_1+n_2=n+1$. Other relevant actions on $\mathbb{S}^n$ satisfying \eqref{OrbitsPositiveDimensional} include the action of $\Gamma=\mathbb{S}^1$ when $n$ is odd, as well as the more intricate examples constructed by Kristaly in \cite{Kristaly2009}. Further important examples for which the symmetry condition \eqref{OrbitsPositiveDimensional} can be realized are provided by the complex projective space $\mathbb{C}\mathbb{P}^n$, equipped with the Fubini--Study metric, and the quaternionic projective space $\mathbb{H}\mathbb{P}^n$, endowed with its standard symmetric metric. In these cases, one may consider, respectively, actions induced by products of unitary groups $\Gamma=U(n_1)\times U(n_2)$ with $n_1+n_2=2n+2$, and by products of symplectic groups $\Gamma=Sp(n_1)\times Sp(n_2)$ with $n_1+n_2=4n+4$. These examples, together with several other less elementary constructions satisfying the same symmetry requirement, are discussed in Section 7 of \cite{FernandezPalmasTorresOrozco2024}. These examples therefore place Theorem \ref{Theorem:Main} in a broad geometric setting and lead to the following immediate corollary, which extends the main result of \cite{BartschLiuWeth2004} to a broader class of Einstein manifolds and group actions, and generalizes Corollary 1.4 of \cite{FernandezPalmasTorresOrozco2024} and Corollary 1.4 of \cite{JulioBatallaPetean2025} to higher-order operators.

\begin{corollary}
Let $\mu>0$, $m\in\mathbb{N}$, $(M,g)$ be a closed Einstein manifold of dimension $n> 2m$ with positive scalar curvature, and $\Gamma$ be a closed subgroup of $\text{Isom}(M,g)$. If $\Gamma$ satisfies the symmetry condition \eqref{OrbitsPositiveDimensional}, then problem \eqref{MainEquation:Q-curvatureEquation} admits an infinite number of $\Gamma$-invariant solutions with increasing energies. One of them is positive,  while the rest change sign.
\end{corollary}

\medskip

As mentioned above, a key feature of Einstein manifolds that allows the operator $\mathscr{P}_g^m$ to satisfy coercivity and property \eqref{Property:PPP} is the existence of a useful factorization into second-order conformally covariant operators. More recently, Case and Malchiodi \cite{CaseMalchiodi2024} obtained an analogous factorization for $\mathscr{P}_g^m$ on special Einstein products, namely, Riemannian products of Einstein manifolds $(N^\ell\times L^{n-\ell},h\oplus k)$ satisfying $\text{Ric}_h=-(\ell-1)\lambda h$ and $\text{Ric}_k=(n-\ell-1)\lambda k$ for some constant $\lambda>0$. In this setting, the authors proved the existence of a constant $D(m,\ell)$, depending only on the order $m$ of the operator and on $\ell$, such that, if $n\geq D(m,\ell)$ and $(N^\ell\times L^{n-\ell},h\oplus k)$ is a closed special Einstein product, then the operator satisfies property \eqref{Property:PPP}. The authors also showed that
\begin{equation}\label{ConstantD(m,l)}
D(m,\ell)\geq 2m+2\ell -1
\end{equation}
and that this is an equality when $m=1,2$, i.e., it is explicit in the case of the conformal laplacian and in the case of the Paneitz-Branson operator, respectively. If the dimension of the closed special Einstein product is large enough, our main result  applies, yielding the following consequence, which, to the best of our knowledge, provides the first existence result for sign-changing solutions to problem \eqref{MainEquation:Q-curvatureEquation} on non-Einstein manifolds.

\begin{corollary}\label{Corllary:SpecialEinsteinProducts}
Let $\mu>0$ and $m,\ell,n\in\mathbb{N}$. Consider a closed special Einstein product $(N^\ell\times L^{n-\ell},h\oplus k)$, and let $\Gamma_N\subset\operatorname{Isom}(N,h)$ and $\Gamma_L\subset\operatorname{Isom}(L,k)$ be closed subgroups. If $n\geq D(m,\ell)$ and the actions of $\Gamma_N$ and $\Gamma_L$ satisfy
\begin{equation}\label{OrbitsSpecialEinsteinProduct}
1\leq \dim\Gamma_N x + \dim\Gamma_Ky< n,\qquad\text{ for all }\; x\in N \text{ and all }y\in L,
\end{equation}
then, the problem \eqref{MainEquation:Q-curvatureEquation} admits an infinite number of $\Gamma_N\times\Gamma_L$-invariant solutions with increasing energies. One of them is positive, and the rest change sign.
\end{corollary}

We next provide explicit examples of special Einstein products where the hypotheses of Corollary \ref{Corllary:SpecialEinsteinProducts} hold true. In what follows, $(\mathbb{S}^\ell,g_{\ell})$ denotes the round sphere of dimension $\ell\geq 1$, $(\mathbb{H}^\ell,\hat{g}_\ell)$ denotes the $\ell$-dimensional hyperbolic space endowed with its canonical metric, and recall that $\text{Isom}(\mathbb{S}^\ell,g_\ell)=O(\ell+1)$, the orthogonal group in dimension $\ell+1$.

\begin{example} Let $m\geq 2$ be a positive integer. Then $(\mathbb{S}^1\times \mathbb{S}^{n-1},g_1\oplus g_{n-1})$ is a closed special Einstein product, and if we consider any closed subgroup $\Gamma_{\mathbb{S}^{n-1}}$ of $O(n)$ with positive dimensional orbits and $\Gamma_{\mathbb{S}^1}=\{Id\}$, then the hypotheses in Corollary \ref{Corllary:SpecialEinsteinProducts} are satisfied for any $n\geq D(m,1)$. An interesting case is given when $m=3$, because, as remarked in \cite[Remark 4.2]{CaseMalchiodi2024}, $D(3,1)\leq 48$ and, hence, equation \eqref{MainEquation:Q-curvatureEquation} for the sixth order operator $\mathscr{P}_{g_1\oplus g_{47}}^3$ on $\mathbb{S}^1\times \mathbb{S}^{47}$, admits an infinite number of $\{Id\}\times O(48)$-invariant sign-changing solutions.
\end{example}

\begin{example}
For any $\ell\geq 2$, a large class of special Einstein products is obtained by taking the product of a closed hyperbolic $\ell$-manifold and the round sphere $\mathbb{S}^{n-\ell}$ (see \cite{CaseMalchiodi2024} and \cite[Example 1.167]{Besse1987}). A natural way to construct closed $\ell$-dimensional hyperbolic manifolds is by considering the quotient $\mathbb{H}^\ell/\Lambda$ by a torsion free and uniform (or cocompact) lattice $\Lambda$ of $\text{Isom}(\mathbb{H}^\ell,\hat{g}_\ell)$, see \cite[Section 4.14]{Kapovich2001}. Explicit and interesting examples for any dimension $\ell\geq 2$ are given, for instance, in \cite{BelolipetskyThomson2011} and in \cite{LongReid2005}. Hence, for any fixed integers $m,\ell\geq 2$, if $n\geq D(m,\ell)$ and if $(N_h^\ell,h)$ denotes a closed hyperbolic manifold, then the special Einstein product $(N^\ell\times \mathbb{S}^{n-\ell}, h\times g_{n-\ell})$, together with the action of $\{Id\}\times \Gamma_{\mathbb{S}^{n-\ell}}$, where $\Gamma_{\mathbb{S}^{n-\ell}}$ is any closed subgroup of $O(n-\ell+1)$ with positive dimensional orbits, satisfy the hypotheses of Corollary \ref{Corllary:SpecialEinsteinProducts} and equation \eqref{MainEquation:Q-curvatureEquation} admits an infinite number of sign-changing solutions. Since inequality \eqref{ConstantD(m,l)} is an equality when $m=2$, taking any $\ell\geq 2$, $n=3+2\ell$ and the special Einstein product $(N^\ell\times \mathbb{S}^{3+\ell})$, we have that the Paneitz-Branson equation, i.e., equation \eqref{MainEquation:Q-curvatureEquation} with $m=2$, admits an infinite number of sign-changing solutions. We want to stress that, as far as can be determined from the current literature, this is the first known example of sign-changing solutions to the Paneitz-Branson $Q$-curvature equation on non-Einstein manifolds.
\end{example}

\bigskip

This paper is organized as follows. In Section \ref{Section:VariationalSetting} we provide the functional setting related to $GJMS$-operators. In Section \ref{Section:VariationalSettingSymmetries}, we introduce the symmetric variational setting. Later, in Section \ref{Section:MultiplicitySIgnChangingSolutions}, we introduce the dual cone method developed in \cite{Weth2006} to prove the existence of sign-changing solutions to the problem \eqref{MainEquation:Q-curvatureEquation}. Finally, in Section \ref{Section:EinsteinProducts}, we give conditions to guarantee that the $GJMS$-operators defined on closed special Einstein products are coercive, and we prove Corollary \ref{Corllary:SpecialEinsteinProducts}.

\section{Variational Setting}\label{Section:VariationalSetting}

Let $m\in\mathbb{N}$ and consider a closed Riemannian manifold $(M,g)$ of dimension $n>2m$.

In what follows, we denote the Lebesgue spaces associated with $(M,g)$ by $L_g^s(M), s\geq 1$, with the norm 
\[
\Vert u\Vert_s:=\left( \int_M \vert u\vert^s \right)^{1/s}.
\]
We define the Sobolev space $H_g^m(M)$ as the completion of $C^\infty(M)$ with respect to the norm
\begin{equation*}\label{Equation:UsualSobolevNorm}
\|v \|_{H_g^{m}(M)}=\left[ \sum_{j=0}^{m}  \int_{M}   |\nabla^{j} v |_g^{2} dV_g  \right]^{1/2}, 
\end{equation*}
endowed with the inner product
\begin{equation*}\label{Equation:UsualSobolevInnerProduct}
\langle u,v \rangle_{H_g^{m}(M)}:= \sum_{j=0}^{m}  \int_{M}   \langle\nabla^{j} u , \nabla^{j}v\rangle_g\; dV_g,\quad u,v\in C^\infty(M) , 
\end{equation*}
where $\nabla^j u$ denotes the $j$-th covariant derivative of $u$ and 
\[
\langle \nabla^j u,\nabla^j v\rangle_g:=g^{i_1k_1}\cdots g^{i_j,k_j}(\nabla^j u)_{i_1\cdots i_j}(\nabla^j v)_{k_1\cdots k_j},
\]
denotes the Riemannian metric extended to tensor fields.

Through this section, we suppose that the $GJMS$-operator $\mathscr{P}_g$ associated with $(M,g)$ is coercive and satisfies  \eqref{Property:PPP}. The self-adjoint property of $P_g$ implies that the bilinear form
\[
\langle u,v \rangle_\ast := \int_M u\mathscr{P}_g^m v\; dV_g,\qquad u,v\in C^\infty(M),
\]
extends to a continuous symmetric bilinear form on $H_g^m(M)\times H_g^m(M)$ (see \cite[Proposition 1]{Robert2011}). The coercivity of the operator yields that this bilinear form is a well defined inner product in $H_g^m(M)$ and the norm
\[
\Vert u\Vert_\ast := \sqrt{\langle u, u\rangle_\ast} = \sqrt{\int_M u \mathscr{P}_g^m u\; dV_g},\quad u\in H_g^m(M),
\]
is equivalent to the norm $\Vert\cdot\Vert_{H_g^m(M)}$ (see \cite[Proposition 2]{Robert2011}).

\begin{remark}\label{Remark:Renormalized}
Notice that $u$ is a solution to \eqref{MainEquation:Q-curvatureEquation} if and only if
\[
v=\mu^{\frac{1}{2-2_{m,n}^\ast}}u 
\] 
is a solution to the renormalized problem
\begin{equation}\label{RenormalizedQ-curvatureEquation}
\mathscr{P}_g^m(v) = \vert v\vert^{2_{m,n}^\ast-2}v\qquad \text{ on } M
\end{equation}

\end{remark}

In light of the preceding remark, we can restrict our attention to equation \eqref{RenormalizedQ-curvatureEquation}. 
Define the energy functional $\Phi: H_g^m(M)\to \mathbb{R}$
\[
\Phi(u):= \frac{1}{2}\Vert u\Vert_\ast^2 - \Psi(u), \quad \Psi(u):= \frac{1}{2_{m,n}^\ast}\int_M \vert u\vert^{2_{m,n}^\ast} \;dV_g.
\]

Since the Sobolev embedding $H_g^m(M)\hookrightarrow L_g^{2_{m,n}^\ast}(M)$ is continuous, $\Psi$ is well defined and of class $C^{2}$. Hence, $\Phi$ is also well defined, of class $C^{2}$ and with derivative given by
\[
\Phi'(u)v = \langle u,v\rangle_\ast - \Psi'(u)v = \langle u,v\rangle_\ast - \int_M \vert u \vert^{2_{m,n}^\ast-2}uv\; dV_g.
\]
Thus, critical points of $\Phi$ correspond to the solutions of \eqref{RenormalizedQ-curvatureEquation}. 
Also notice that
\begin{equation}\label{Eq:PsiEven}
\Psi(-u) = \Psi(u)\qquad \text{and}\qquad \Phi(-u)=\Phi(u).
\end{equation}

The dual cone method locates nontrivial critical points of $\Phi$ applying dynamical system arguments to the negative gradient flow of $\Phi$. To do so, define the nonlinear operator $A:H_g^m(M)\to H_g^m(M)$ by $A(u):=\nabla\Psi(u)$, that is
\begin{equation}\label{Eq:DefinitionA}
\langle A(u),v\rangle_\ast = \langle \nabla\Psi(u),v\rangle_\ast = \Psi'(u)v = \int_M \vert u\vert^{2_{m,n}^\ast-2}uv \; dV_g,\qquad\forall \;u,v\in H_g^m(M).
\end{equation}
Hence, the existence of critical points of $\Phi$ is equivalent to finding fixed points of $A$. 

Observe that $A$ is Lipschitz continuous because $\Psi$ is a $C^{2}$-functional. Moreover, the operator $A$ satisfy the following properties.

\begin{lemma}\label{Lemma:A2}
\begin{enumerate}[label=(\roman*)]
\item \label{Item:A_Odd} $A(-u)=-A(u)$ for every $u\in H_g^m(M)$ and $A(0)=0$.
\item \label{Item:A_BoundedBelow} There exist constants $C_0>0$ and $\eta>2$  such that
\begin{align*}
\langle A(u), u\rangle_\ast&\geq \eta\Psi(u) - C_0, & \forall\; u\in H_g^m(M)
\end{align*}
\item \label{Item:A_Grow}  There exist constants $q^\ast>0$, $p>1$, and $q_\ast\in(0,1)$ such that
\begin{align*}
\vert \langle A(u), v\rangle_\ast\vert &\leq (q_\ast\Vert u\Vert_\ast + q^\ast\Vert u\Vert_\ast^{p})\Vert v\Vert_{\ast}, & \forall\; u,v\in H_g^m(M).
\end{align*}
\end{enumerate}
\end{lemma}

\begin{proof}
Item \ref{Item:A_Odd} follows directly from \eqref{Eq:DefinitionA}: indeed, given $u\in H_g^m(M)$, for any $v\in H_g^m(M)$, we have that
\[
\langle A(-u),v \rangle_\ast = \int_M \vert -u\vert^{2_{m,n}^\ast-2}(-u)v \; dV_g = -\int_M \vert u\vert^{2_{m,n}^\ast-2}uv \; dV_g = \langle -A(u) , v \rangle_\ast,
\]
where we conclude that $A(-u)=-A(u)$ and, consequently, $A(0)=0$.

We now turn to the proof of \ref{Item:A_BoundedBelow}. Again, by \eqref{Eq:DefinitionA}, for any $C_0>0$, we have that
\[
\langle A(u),u \rangle_\ast = \int_M \vert u\vert^{2_{m,n}^\ast-2}u^2\; dV_g = 2_{m,n}^\ast\Psi(u) > 2_{m,n}^\ast\Psi(u) - C_0.
\]
Thus, the inequality follows by taking $\eta:=2_{m,n}^\ast>2$.

Finally, we prove \ref{Item:A_Grow}. From the Sobolev embedding $H_g^m(M)\hookrightarrow L_g^{2_{m,n}^\ast}(M)$, there exists a constant $C>0$ such that
\begin{equation}\label{Eq:CriticalSobolevEmbeding}
\Vert u \Vert_{2_{m,n}^\ast}\leq C \Vert u\Vert_\ast,\qquad\forall\; u\in H_g^m(M).
\end{equation}
Therefore, from \eqref{Eq:DefinitionA}, together with Hölder inequality for $s:=2_{m,n}^\ast$ and $\frac{s}{s-1}$, and \eqref{Eq:CriticalSobolevEmbeding}, taking $q_\ast\in(0,1)$ arbitrary, for every $u,v\in H_g^m(M)$ we have that
\begin{align*}
\left\vert  \langle A(u) , v\rangle_\ast \right\vert &\leq \int_M \vert u\vert^{2_{m,n}^\ast-1}\vert v\vert \; dV_g\\
&\leq \left(\int_M \vert u\vert^{2_{m,n}^{\ast}}\; dV_g\right)^{\frac{2_{m,n}^\ast-1}{2_{m,n}^\ast}} \left(\int_M \vert v\vert^{2_{m,n}^{\ast}}\; dV_g\right)^{\frac{1}{2_{m,n}^\ast}}\\
&= \Vert u\Vert_{2_{m,n}^\ast}^{2_{m,n}^\ast-1}\Vert v\Vert_{2_{m,n}^\ast}\\
&\leq C^{2_{m,n}^\ast}\Vert u\Vert_\ast^{2_{m,n}^\ast-1}\Vert v\Vert_\ast\\
&\leq \left(  q_\ast \Vert u\Vert_\ast + C^{2_{m,n}^\ast}\Vert u\Vert_\ast^{2_{m,n}^\ast-1} \right)\Vert v\Vert_\ast.
\end{align*}
As $2_{m,n}^\ast>2$, taking $p=2_{m,n}^\ast-1>1$ and $q^\ast:=C^{2_{m,n}^\ast}>0$, we conclude the proof of item \ref{Item:A_Grow}.
\end{proof}

A crucial requirement for the existence of fixed points is for the operator $A$ to be compact. However, the presence of the critical Sobolev nonlinearity prevents this from happening. In the next section, we impose symmetry assumptions in order to recover the compactness of the operator $A$ and verify that it satisfies the hypotheses of Theorem 2.7 in \cite{Weth2006}, proving our main result.


\section{Variational Setting with symmetries} \label{Section:VariationalSettingSymmetries}

Throughout this section, we fix a closed subgroup $\Gamma$ of $\text{Isom}(M,g)$ satisfying \eqref{OrbitsPositiveDimensional}. The group $\Gamma$  acts on $H_{g}^{m}(M)$ in the usual way:
\begin{align*}
    \Gamma \times H_{g}^{m}(M) \to H_{g}^{m}(M), \quad (\gamma, u)\mapsto \gamma u:=u \circ \gamma^{-1}.
\end{align*}

The space of fixed points under this action coincides with the subspace of $\Gamma$-invariant functions 
\begin{align*}
H_{g}^{m}(M)^{\Gamma}&:=  \left\{ u \in H^{m}_{g}(M)\;:\, u \text{ is }\Gamma\text{-invariant}  \right\}\\
&= \left\{ u \in H^{m}_{g}(M)\;:\, \gamma u=u \text{ for every }\gamma\in \Gamma  \right\}
\end{align*}
which is closed in $H_g^m(M)$. 
Moreover, if $C^{\infty}(M)^{\Gamma}$ denotes the space of smooth $\Gamma$-invariant functions, $H_{g}^{m}(M)^{\Gamma}$ coincides with the closure of this space under the Sobolev norm $\Vert \cdot\Vert_{H_g^m(M)}$. As $M$ is compact, it follows from the Myers-Steenrod Theorem that $\text{Isom}(M,g)$ is compact (see, for instance, Section II.1 in \cite{Kobayashi1995}). Since $\Gamma$ is closed in $\text{Isom}(M,g)$, $\Gamma$ is compact and admits a Haar measure. This fact and the inequality $\dim\Gamma x\leq n-1$ for all $x\in M$, given in \eqref{OrbitsPositiveDimensional}, imply the existence of $\Gamma$-invariant partitions of unity, see \cite{Palais1961}, and hence, the spaces $C^{\infty}(M)^{\Gamma}$ and $H_{g}^{m}(M)^{\Gamma}$ are infinite-dimensional.   

For any $\gamma\in\text{Isom}(M,g)$, the action of $\Gamma$ on $H_g^m(M)$ induces linear maps
\begin{align*}
\gamma:  H_{g}^{m}(M) \to H_{g}^{m}(M), \quad u \to \gamma u:=u \circ \gamma^{-1}.
\end{align*}
Clearly, each map is a linear isomorphism. We verify that they are isometric isomorphisms.

\begin{lemma}\label{Lemma:P_gCommutes}
 For every  $\gamma\in \text{Isom}(M,g)$ and every $u \in C^{\infty}(M)$,
 \begin{align*}\label{Equation:P_gCommutes}
     \mathscr{P}_{g}^m(u \circ \gamma)=(\mathscr{P}_{g}^mu)\circ \gamma
 \end{align*}
 In particular, 
 \begin{align*}
\gamma:  (H_{g}^{m}(M),\Vert\cdot\Vert_{\ast}) \to (H_{g}^{m}(M),\Vert\cdot\Vert_{\ast})
\end{align*}
is a linear isometry, i.e.,
\[
\langle u\circ\gamma,v\circ\gamma\rangle_{\ast}= \langle u,v\rangle_{\ast},\qquad\forall\;u,v\in H_g^m(M)\text{ and }\forall\; \gamma\in\Gamma.
\]
\end{lemma}

\begin{proof}
	Let $u\in\mathcal{C}^\infty(M)$ and $\gamma\in\text{Isom}(M,g)$. Then $\gamma:M\rightarrow M$ is a diffeomorphism and, by the naturality of $\mathscr{P}_{g}^m$, property \hyperref[Propoerty:Pg naturality]{($P_2$)} in the introduction, we obtain for any isometry $\gamma\in \text{Isom}_g(M,g)$ that
	\[
	\mathscr{P}_{g}^m(u\circ\gamma) = (\mathscr{P}_{g}^m\circ\gamma^\ast)(u)=\mathscr{P}_{\gamma^\ast g}^m\circ\gamma^\ast(u)=(\gamma^\ast\circ \mathscr{P}_{g}^m)(u)=\mathscr{P}_{g}^m(u)\circ\gamma,
	\]
	where $\gamma^\ast g$ denotes the pullback metric.
	
	Now, to see that $\gamma$ induces a linear isometry, recall that if $\gamma\in\text{Isom}(M,g)$, then, by the Change of Variables Theorem, 
	\begin{equation*}\label{Eq:Integral invariance isometries}
		\int_M u\circ\gamma\; dV_g = \int_M u\;  dV_g,
	\end{equation*}
	(Cf. \cite[Théorème 4.1.2]{Hebey1997}).
	Then, for every $u,v\in\mathcal{C}^\infty(M)$,
	\[
	\langle  u\circ\gamma^{-1},v\circ\gamma^{-1}\rangle_\ast=\int_M (v\circ\gamma^{-1})\mathscr{P}_{g}^m(u\circ\gamma^{-1}) dV_g =\int_M (v\circ\gamma^{-1})\mathscr{P}_{g}^m(u)\circ\gamma^{-1} \,dV_g = \int_M v\mathscr{P}_{g}^mu \,dV_g =\langle  u, v\rangle_\ast.
	\]
Density of $\mathcal{C}^\infty(M)$ in $H_g^m(M)$ yields the result.
\end{proof}

\begin{lemma}\label{Lemma:GammaInvariant}
The functionals $\Psi$ and $\Phi$ are $\Gamma$-invariant. Moreover,
\begin{equation}\label{RestrictionWellDefined}
\nabla\Psi(u)\in H_g^m(M)^\Gamma,\qquad\forall\; u\in H_g^m(M)^\Gamma
\end{equation}
\end{lemma}

\begin{proof}

By the Change of Variable Theorem, we  have that
\[
\Psi(u\circ \gamma) = \frac{1}{2_{m,n}^\ast}\int_M \vert u\circ \gamma\vert^{2_{m,n}^\ast}\; dV_g = \frac{1}{2_{m,n}^\ast}\int_M \vert u\vert^{2_{m,n}^\ast}\; dV_g = \Psi(u),
\]
for any $u\in H_g^m(M)$ and any $\gamma\in\Gamma$, concluding that $\Psi$ is $\Gamma$-invariant. 
By Lemma \ref{Lemma:P_gCommutes}, for any $u\in H_g^m(M)$ and any $\gamma\in \Gamma$ we have that $\Vert u\circ\gamma\Vert_\ast=\Vert u\Vert_\ast$.
From this, we conclude that $\Phi(u\circ\gamma)=\Phi(u)$ for any $u\in H_g^m(M)$ and any $\gamma\in\Gamma$; thus, $\Phi$ is also $\Gamma$-invariant.

Now we prove the second assertion. As $\Psi$ is $\Gamma$-invariant, we have that $\Psi=\Psi\circ\gamma$ for any $\gamma\in\Gamma$. Since $\gamma:H_g^m(M)\to H_g^m(M)$ is a linear isometry, then $\gamma'(u)=\gamma$ for any $u\in H_g^m(M)$. Hence, by the chain rule,
\[
\Psi'(u)=(\Psi\circ\gamma)'(u) =\Psi'(\gamma u)\circ\gamma'(u) = \Psi'(\gamma u)\circ\gamma.
\]
Therefore, by Lemma \ref{Lemma:P_gCommutes}, given $u\in H_g^m(M)$ and $\gamma\in\Gamma$,
\begin{align*}
\langle \nabla\Psi(u), v \rangle_\ast = \Psi'(u) v = \Psi'(\gamma u)\gamma v = \langle \nabla \Psi(\gamma u), \gamma v\rangle_\ast= \langle \gamma^{-1} \nabla \Psi(\gamma u),  v\rangle_\ast,\qquad \forall\; v\in H_g^m(M),
\end{align*}
where we obtain that $\nabla\Psi(\gamma u)=\gamma\Psi(u)$. From this, if $u\in H_g^m(M)^\Gamma$, then for any $\gamma\in\Gamma$, $\gamma u = u\circ\gamma^{-1}=u$, concluding that $\gamma\nabla\Psi(u)=\nabla\Psi(u)$ and that $\nabla\Psi(u)\in H_g^m(M)^\Gamma$, as we wanted to prove.

\end{proof}

As a consequence of the last lemma, by the Principle of Symmetric Criticality \cite{Palais1979}, the critical points of $\Phi$ restricted to $H_g^m(M)^\Gamma$ correspond to the $\Gamma$-invariant solutions to the problem \eqref{MainEquation:Q-curvatureEquation}. Hence, we can restrict ourselves to finding critical points of the restriction $\Phi: H_g^m(M)^\Gamma\to\mathbb{R}$. By \eqref{RestrictionWellDefined}, we have that the restriction of $A$ to $H_g^m(M)^\Gamma$ satisfies that
\[
A: H_g^m(M)^\Gamma\to H_g^m(M)^\Gamma
\]
and the critical points of $\Phi$ restricted to $H_g^m(M)^\Gamma$ correspond to the fixed points of $A$ in $H_g^m(M)^\Gamma$. To prove that the restricted operator recovers the compactness, we need the following classical result (see \cite{HebeyVaugon1997} for the case $m=1$ and, for instance, \cite[Lemma 2.2]{FernandezPalmasTorresOrozco2024} for $m\geq 2$).

\begin{lemma}\label{Lemma:HebeyVaugon}
    If $\kappa:=\min \{ \dim\; \Gamma x :\,\, x \in M  \}$, then 
    \begin{align*}
        H_{g}^{m}(M)^{\Gamma} \hookrightarrow L_{g}^{r}(M)
    \end{align*}
    is continuous for every $2<r \leq 2^{\ast}_{m,n-\kappa}:= \frac{2(n-\kappa)}{(n-\kappa)-2m}$ and  compact if $2<r < 2^{\ast}_{m,n-\kappa}$.
\end{lemma}

We next prove that $A$ is compact.

\begin{proposition}\label{Proposition:A1}
The operator $A:H_g^m(M)^\Gamma\to H_g^m(M)^\Gamma$ is compact.
\end{proposition}

\begin{proof}
Let $(u_k)$ be a bounded sequence in $H_g^m(M)^\Gamma$. Since $\dim \Gamma x\geq 1$ by \eqref{OrbitsPositiveDimensional}, we have that $\kappa=\min_{x\in M}\dim\Gamma x\geq 1$ and $2_{m,n}^\ast < 2_{m,n-\kappa}^\ast$. Hence, by Lemma \ref{Lemma:HebeyVaugon}, the Sobolev embedding $H_g^m(M)^\Gamma\hookrightarrow L^{s^\ast}_g(M)$ is compact for $s=2$ and $s=2_{m,n}^\ast$. Then, up to a subsequence, which we denote the same, there exists $u_0\in H_g^m(M)^\Gamma$ such that
\begin{align*}
u_k\rightharpoonup u_0 &\qquad \text{weakly in }H_g^m(M)^\Gamma\\
u_k\to u_0 &\qquad \text{strongly in }L_g^{2_{m,n}^\ast}(M).
\end{align*}

We next show that $A(u_k)\to A(u_0)$ in $H_g^m(M)$ in a subsequence. First, as $2_{m,n}^\ast>2$, there exists $C>0$ such that
\[
\vert \vert a\vert^{2_{m,n}^\ast-2}a - \vert b\vert^{2_{m,n}^\ast-2}b\vert \leq C (\vert a\vert + \vert b\vert)^{2_{m,n}^\ast-2}\vert a- b\vert,\qquad\forall\;a,b\in\mathbb{R},
\]
see \cite[Lemma 1]{DincaJebeleanMawhin2001}. Then, as $u_k,u_0\in L^{2_{m,n}^\ast}_g(M)$,  $\vert u_k\vert^{2_{m,n}^\ast-2}, \vert u_0\vert^{2_{m,n}^\ast-2}\in L_g^{\frac{2_{m,n}^\ast}{2_{m,n}^\ast-2}}(M)$, and as  $\frac{2_{m,n}^\ast - 2}{2_{m,n}^\ast}+\frac{1}{2_{m,n}^\ast}=\frac{2_{m,n}^\ast-1}{2_{m,n}^\ast}$, the generalized Hölder inequality yields that
\begin{align*}
\left[\int_{M} \left\vert  \vert u_k\vert^{2_{m,n}^\ast-2}u_k - \vert u_0\vert^{2_{m,n}^\ast-2}u_0 \right\vert^{\frac{2_{m,n}^\ast}{2_{m,n}^\ast-1}}\; dV_g \right]^{\frac{2_{m,n}^\ast -1}{2_{m,n}^\ast}}
&\leq \left[\int_M \left[ (\vert u_k\vert + \vert u_0\vert)^{2_{m,n}^\ast-2} \vert u_k - u_0\vert \right]^{\frac{2_{m,n}^\ast}{2_{m,n}^\ast-1}}\;dV_g\right]^{\frac{2_{m,n}^\ast -1}{2_{m,n}^\ast}}\\
&\leq \left(\int_M (\vert u_k\vert + \vert u_0\vert)^{2_{m,n}^\ast} \; dV_g\right)^{\frac{2_{m,n}^\ast-2}{2_{m,n}^\ast}}\left( \int_M \vert u_k- u_0\vert^{2_{m,n}^\ast} \right)^{\frac{1}{2_{m,n}^\ast}}\\
&\leq C\Vert u_k - u_0\Vert_{2_{m,n}^\ast},
\end{align*}
where $C>0$ is a positive constant; its existence is guaranteed because $u_k\to u_0$ in $L_g^{2_{m,n}^\ast}(M)$ and is therefore bounded. From this inequality and Hölder inequality for $\frac{2_{m,n}^\ast}{2_{m,n}^\ast-1}$ and $2_{m,n}^\ast$, for every $v\in H_g^m(M)^\Gamma$ we have that
\begin{align*}
\vert \langle A(u_k)-A(u_0),v\rangle_\ast\vert 
&\leq \int_M \left\vert  \vert u_k\vert^{2_{m,n}^\ast-2}u_k - \vert u_0\vert^{2_{m,n}^\ast-2}u_0 \right\vert\vert v\vert \;dV_g\\
&\leq \left( \int_M\left\vert  \vert u_k\vert^{2_{m,n}^\ast-2}u_k - \vert u_0\vert^{2_{m,n}^\ast-2}u_0 \right \vert^{\frac{2_{m,n}^\ast}{2_{m,n}^\ast-1}} \; dV_g \right)^{\frac{2_{m,n}^\ast-1}{2_{m,n}^\ast}}\left(\int_M \vert v\vert^{2_{m,n}^\ast}\; dV_g\right)^{\frac{1}{2_{m,n}^\ast}}\\
&\leq C \Vert u_k - u_0\Vert_{2_{m,n}^\ast}\Vert v\Vert_{2_{m,n}^\ast}\to 0,
\end{align*}
because $u_k\to u_0$ in $L_g^{2_{m,n}^\ast}(M)$ in a subsequence. Hence, $A(u_k)\to A(u_0)$ in a subsequence and the operator $A$ is compact.

\end{proof}

\begin{corollary}
The functional $\Phi$ satisfies the Palais-Smale condition in $H_g^m(M)^\Gamma$, i.e., every  sequence $(u_k)$ in  $H_g^m(M)^\Gamma$ such that $(\Phi(u_k))$ is bounded in $\mathbb{R}$ and $\nabla\Phi(u_k)\to0$ as $n\to\infty$, then $(u_k)$ has a convergent subsequence.
\end{corollary}

\begin{proof}
This is a straightforward consequence of Proposition \ref{Proposition:A1} and Lemma \ref{Lemma:A2}. See Lemma 2.1 in \cite{Weth2006} for the details.
\end{proof}

Another consequence of the compactness lemma is the following behavior at infinity of the functional $\Phi$ restricted to the space $H_g^m(M)^\Gamma.$

\begin{lemma}\label{Lemma:FiniteDimensional}
 If $S\subset H_g^m(M)^\Gamma\smallsetminus\{0\}$ is a compact subset and $\widetilde{S}:=\{tu\;:\; u\in S, t\geq 0\}$, then
\[
\lim_{u\in\widetilde{S},\Vert u\Vert_\ast\to\infty}\Phi(u)=-\infty.
\]
\end{lemma}

\begin{proof}
It suffices to prove that any sequence $(u_k)$ in $\widetilde{S}$ such that $\Vert u_k\Vert_\ast\to\infty$ as $k\to\infty$, has a subsequence, which we denote the same, such that $\Phi(u_k)\to -\infty$. Indeed, if $(u_k)$ is such a subsequence, then, by definition of $\widetilde{S}$, for each $k\in\mathbb{N}$, there exist $v_k\in S$ and $t_k\geq 0$ such that $u_k=t_kv_k$. As $S$ is compact, then the sequence $(v_k)$ is bounded in $H_g^m(M)^\Gamma$ and there is $v_0\in S$ such that $v_k\to v_0$ in $H_g^m(M)$ in a subsequence. Moreover, as the embedding $H_g^m(M)^\Gamma\hookrightarrow L_g^{2_{m,n}^\ast}(M)$ is compact by Lemma \ref{Lemma:HebeyVaugon}, passing again to a subsequence, we have that $v_k\to v_0$ in $L^{2_{m,n}^\ast}_g(M)$. Consequently, up to a subsequence,
\[
\Vert v_k\Vert_\ast^2 =\Vert v_0\Vert_\ast^2 + o(1)\qquad\text{and}\qquad\Vert v_k\Vert_{2_{m,n}^\ast}^{2_{m,n}^\ast}= \Vert v_0\Vert_{2_{m,n}^\ast}^{2_{m,n}^\ast}+ o(1),
\]
where $o(1)\to0$ as $k\to\infty$. On the other hand, as $(v_k)$ is bounded in $H_g^m(M)$ and $\Vert u_k\Vert_\ast\to\infty$, it follows that $t_k\to\infty$. With this, and as $2_{m,n}^\ast>2$, we have that
\begin{align*}
\Phi(u_k) &= \frac{1}{2}\Vert u_k\Vert_\ast^2 - \frac{1}{2_{m,n}^\ast}\Vert u_k\Vert_{2_{m,n}^\ast}^{2_{m,n}^\ast} \\
&= \frac{t_k^2}{2}\Vert v_k\Vert_\ast^2 - \frac{t_k^{2_{m,n}^\ast}}{2_{m,n}^\ast}\Vert v_k\Vert_{2_{m,n}^\ast}^{2_{m,n}^\ast} \\
& = t_k^{2_{m,n}^\ast} \left[\frac{ t_{k}^{2-2_{m,n}^\ast}}{2}(\Vert v_0\Vert^{2}_\ast + o(1)) - \frac{1}{2_{m,n}^\ast}(\Vert v_0\Vert_{2_{m,n}^\ast}^{2_{m,n}^\ast} + o(1)) \right]\\
&= t_k^{2_{m,n}^\ast}\left[ o(1) -\frac{1}{2_{m,n}^\ast}\Vert v_0\Vert_{2_{m,n}^\ast}^{2_{m,n}^\ast} \right]\\
&\to-\infty\qquad \text{as }k\to\infty.
\end{align*}
Hence, $\Phi(u_k)\to-\infty$ in a subsequence, as we wanted to show.
\end{proof}

We end this section with the following existence result.

\begin{lemma}\label{Lemma:ExistenceUniquenessGammaInvariant}
For any $\Gamma$-invariant function $h\in L^2_g(M)$, the linear problem
\begin{equation}\label{LinearProblem}
\mathscr{P}_g^m v = h\quad \text{ on }\ M
\end{equation}
has a unique $\Gamma$-invariant solution $v\in H_g^m(M)^\Gamma.$ Moreover, if $h\in C^\infty(M)^\Gamma$, then $v\in C^\infty(M)^\Gamma$.
\end{lemma}

\begin{proof}
For any $h\in L_g^2(M)$, as $H_g^m(M)^\Gamma$ endowed with the inner product $\langle\cdot,\cdot\rangle_\ast$ is a Hilbert space, and as the linear mapping $T:H_g^m(M)^\Gamma\to\mathbb{R}$ given by $u\mapsto \int_M hu\; dV_g$ is continuous, by Frechet-Riez Representation Theorem, there exists a unique $v\in H_g^m(M)^\Gamma$ such that
\begin{equation}\label{Equation:RestrictionJ}
\langle v,u \rangle_\ast = \int_M hu\; dV_g,\qquad \forall\; u\in H_g^m(M)^\Gamma
\end{equation}
To see that $v$ is a $\Gamma$-invariant solution to the linear problem \eqref{LinearProblem}, observe that  the solutions to this  problem coincide with the critical points of the $C^1$ functional defined as
\[
J: H_g^m(M)\to\mathbb{R},\qquad J(u):=\frac{1}{2}\Vert u\Vert_\ast^2 - \int_M hu \; dV_g,
\]
which is well defined for every $h\in L_g^2(M)$. In particular, when $h$ is $\Gamma$-invariant,  Lemma \ref{Lemma:P_gCommutes} yields that $J$ is $\Gamma$-invariant and identity \eqref{Equation:RestrictionJ} says that $v\in H_g^m(M)^\Gamma$ is a critical point of $J$ restricted to $H_g^m(M)^\Gamma$. Hence, the Priciple of Symmetric Criticality yields that $v$ is a critical point of $J$ in $H_g^m(M)$ and an actual $\Gamma$-invariant (weak) solution to the problem \eqref{LinearProblem}, as we wanted to prove.

The assertion about the smoothness of $v$ given the smoothness of $h$ is a standard regularity argument, because $\mathscr{P}_g^m$ is an elliptic differential operator of order $2m$ defined on a compact manifold, see, for instance \cite[Chapter 5, Theorem 11.1]{Taylor2023}.
\end{proof}

\section{Multiplicity of sign-changing solutions}\label{Section:MultiplicitySIgnChangingSolutions}

We briefly develop the concepts necessary to apply the dual cone method given in \cite{Weth2006}. Consider the cone of positive functions in $H_g^m(M)^\Gamma$,
\[
\mathcal{K}^\Gamma:= \{u\in H_g^m(M)^\Gamma \;:\; u\geq 0 \text{ a.e. in }M\},
\]
which is a closed cone in $H_g^m(M)^\Gamma$, so that $\mathcal{K}^\Gamma=\overline{\mathcal{K}^\Gamma}$ is convex, $\mathbb{R}^+\cdot \mathcal{K}^\Gamma\subset \mathcal{K}^\Gamma$  and $\mathcal{K}^\Gamma\cap(-\mathcal{K}^\Gamma)=\{0\}$. Let $\mathcal{P}: H_g^m(M)^\Gamma\to\mathcal{K}^\Gamma$ be the projection onto $\mathcal{K}^\Gamma$, determined uniquely by
\[
\Vert u - \mathcal{P}(u)\Vert_\ast = \min_{v\in\mathcal{K}^\Gamma} \Vert u -v \Vert_\ast.
\]
Define the dual cone
\[
\left(\mathcal{K}^\Gamma\right)^\ast:= \{u\in H_g^m(M)^\Gamma\;:\; \langle u,v\rangle_\ast\leq 0,\ \forall\; v\in\mathcal{K}^\Gamma\},
\]
and $\mathcal{P}^\ast:= Id - \mathcal{P}$, so that
\[
\mathcal{P}^\ast(u)\in \left(\mathcal{K}^\Gamma\right)^\ast\quad\text{and}\quad \langle \mathcal{P}(u), \mathcal{P}^\ast(u) \rangle_\ast=0,\quad\forall u\in H_g^m(M)^\Gamma.
\]

We need the following preliminary result about $\Gamma$-invariant solutions to linear problems involving $\mathscr{P}_g^m$.

As usual, given $u:M\to \mathbb{R}$, denote by $u^+$ to the positive part of $u$ and by $u^-$ the negative part, so that $u=u^+ + u^-$. We have the following consequence of the  \eqref{Property:PPP}.

\begin{lemma}\label{Lemma:DualConesPositivity}
\begin{enumerate}[label=(\roman*)]
\item \label{Item:DualCone}$\left(\mathcal{K}^\Gamma\right)^\ast\subset -\mathcal{K}^\Gamma.$
\item $\left(-\mathcal{K}^\Gamma\right)^\ast\subset\mathcal{K}^\Gamma.$
\item $\mathcal{P}(u)\geq u^+$ and $\mathcal{P}^\ast(u)\leq u^-$ for every $u\in H_g^m(M)^\Gamma.$
\end{enumerate}
\end{lemma}

\begin{proof}
It is readily seen that the three conditions are equivalent. So, it suffices to prove that \ref{Item:DualCone} holds true. 

Let $u\in(\mathcal{K}^\Gamma)^\ast$ and consider the auxiliary problem
\begin{equation*}
\mathscr{P}_g^m v = h\quad \text{ on }\ M
\end{equation*}
where $h\in C^\infty(M)^\Gamma\cap\mathcal{K}^\Gamma$. As $h$ is smooth, by Lemma \ref{Lemma:ExistenceUniquenessGammaInvariant}, this problem has a unique solution $v\in C^\infty(M)^\Gamma$ and, as $\mathscr{P}_g^mv=h\geq 0$, by  \eqref{Property:PPP}, then also $v\geq 0$. Since $u\in(\mathcal{K}^\Gamma)^\ast$, it follows that
\[
0\geq \langle u,v \rangle_\ast = \int_M u\mathscr{P}_g^m v \; dV_g= \int_M uh  \; dV_g.
\]
Since $h\geq 0$ was an arbitrary function in $C^\infty(M)^\Gamma$, then $u\leq 0$ a.e. in $M$, concluding that $u\in -\mathcal{K}^\Gamma$, as we wanted to prove.
\end{proof}

\begin{proposition}\label{Proposition:A3}
The operator $A:H_g^m(M)^\Gamma\to H_g^m(M)^\Gamma$ satisfies that
\[
\langle A(u),v \rangle_\ast \leq \langle A(\mathcal{P}^\ast u), v \rangle_\ast,\quad\text{for every }u\in H_g^m(M)^\Gamma \text{ and every }v\in\left(\mathcal{K}^\Gamma\right)^\ast.
\]
\end{proposition}

\begin{proof}
Let $u\in H_g^m(M)^\Gamma$ and $v\in \left(\mathcal{K}^\Gamma\right)^\ast$. Then,  $\mathcal{P}^\ast u\in \left(\mathcal{K}^\Gamma\right)^\ast$ and by Lemma \ref{Lemma:DualConesPositivity} we have that $\mathcal{P}^\ast u\leq 0$ and $\mathcal{P}^\ast u\leq u^-\leq u$. As $2_{m,n}^\ast>2$, the function $t\to \vert t\vert^{2_{m,n}^\ast-2}t$ is increasing and, on the one hand, we get that
\begin{align}\label{Eq:InequalityP}
\vert \mathcal{P}^\ast u\vert^{2_{m,n}^\ast-2}\mathcal{P}^\ast u \leq \vert u^-\vert^{2_{m,n}^\ast-2}u^{-}\leq \vert u\vert^{2_{m,n}^\ast-2}u\quad\text{ a.e. in }\ M
\end{align}
On the other hand, since $v\in \left(\mathcal{K}^\Gamma\right)^\ast$, by Lemma \ref{Lemma:DualConesPositivity} it follows that $v\in -\mathcal{K}^\Gamma$ and $v\leq 0$. Hence, from \eqref{Eq:InequalityP} we get that
\[
\vert u\vert^{2_{m,n}^\ast-2}uv\leq \vert u^{-}\vert^{2_{m,n}^\ast-2}u^{-} v\leq \vert \mathcal{P}^\ast u\vert^{2_{m,n}^\ast-2}\mathcal{P}^\ast u v\qquad\text{a.e. in }\ M.
\]
This implies that
\begin{align*}
\langle A(u), v \rangle_\ast = \int_M \vert u\vert^{2_{m,n}^\ast-2}uv\; dV_g
\leq \int_M \vert u^-\vert^{2_{m,n}^\ast-2}u^-v\; dV_g
\leq \int_M \vert \mathcal{P}^\ast u\vert^{2_{m,n}^\ast-2}\mathcal{P}^\ast uv\; dV_g
= \langle A(\mathcal{P}^\ast u), v \rangle_\ast, 
\end{align*}
and the lemma follows.
\end{proof}

We can conclude the existence of a $\Gamma$-invariant positive solution and an infinite number of sign-changing $\Gamma$-invariant solutions to the problem \eqref{RenormalizedQ-curvatureEquation}.


\begin{proof}[Proof of Theorem \ref{Theorem:Main}]

By Remark \ref{Remark:Renormalized}, it suffices to show the result for the reonormalized equation \eqref{RenormalizedQ-curvatureEquation}, and as $\Phi$ is $\Gamma$-invariant, as remarked after Lemma \ref{Lemma:GammaInvariant}, it suffices to show the existence of critical points of $\Phi$ restricted to $H_g^m(M)^\Gamma.$

To prove the existence of sign-changing solutions using the dual cone method, we next verify that all the hypotheses in Theorem 2.7 in \cite{Weth2006} are satisfied. First by Propositions \ref{Proposition:A1} and \ref{Proposition:A3} and by Lemma \ref{Lemma:A2}, the assumptions $(A_1)-(A_3)$ in \cite{Weth2006} hold true. Moreover, by \eqref{Eq:PsiEven}, the functional $\Psi$ is even and by item \ref{Item:A_Odd} in Lemma \ref{Lemma:A2}, the operator $A$ is odd. As remarked at the beginning of this section, the space $H_g^m(M)^\Gamma$ is infinite-dimensional and, for any $k\in\mathbb{N}$, let $V_k$ be a subspace of $H_g^m(M)^\Gamma$ of dimension $k$. Taking $\widetilde{S}=V_k$ and $S=S_{V_k}^1(0)$, the sphere of radius $1$ in $V_k$, in Lemma  \ref{Lemma:FiniteDimensional}, we have that
\[
\limsup_{u\in V_k,\Vert u\Vert_\ast\to\infty}\Phi (u) = \lim_{u\in V_k,\Vert u\Vert_\ast\to\infty}\Phi (u)=-\infty.
\]
In this way, Theorem 2.7 in \cite{Weth2006} can be applied and $\Phi$ possesses, at least, $k-1$ pairs of of nontrivial critical points located in $H_g^m(M)^\Gamma\smallsetminus\left(\mathcal{K}^\Gamma\cup-\mathcal{K}^\Gamma\right)$, i.e., they are sign-changing, and having increasing energies. As $k\in\mathbb{N}$ was arbitrary, then there is an infinite number of pairs of sign-changing critical points of $\Phi$, yielding the existence of an infinite number of pairs of sign-changing $\Gamma$-invariant solutions to the problem  \eqref{RenormalizedQ-curvatureEquation}, with increasing energies.

Finally, the existence of a positive $\Gamma$-invariant solution follows directly from Proposition 2.4 in \cite{Weth2006}, because
for every $u_0\in\mathcal{K}^\Gamma\smallsetminus\{0\}$ we have that
\[
\Phi(tu_0)=\frac{t^2}{2}\Vert u_0\Vert_\ast - \frac{t^{2_{m,n}^\ast}}{2_{m,n}^\ast}\int_M \vert u_0\vert^{2_{m,n}^\ast}\; dV_g\to-\infty,\qquad t\to\infty,
\]
yielding that $\Phi$ has a nontrivial critical point in $\mathcal{K}^\Gamma$ and, thus, a positive $\Gamma$-invariant solution to \eqref{RenormalizedQ-curvatureEquation}.
\end{proof}

\begin{remark}
The existence of a least energy positive $\Gamma$-invariant solution to problem \eqref{MainEquation:Q-curvatureEquation} was proven by Robert  in \cite{Robert2011}. The dual cone method only states the existence of a positive $\Gamma$-invariant solution, and it is not clear to us if the solution coincides with the minimizer obtained by Robert or not. It would be interesting to compare both solutions.
\end{remark}

\section{Coercivity of the GJMS-operator on special Einstein products}\label{Section:EinsteinProducts}

In this section, we prove Corollary \ref{Corllary:SpecialEinsteinProducts}. To do this, we need to provide conditions to guaranty that the $GJMS$-operator on an special Einstein product is positive. Recall from the introduction the definition of a closed special Einstein product $(N^\ell\times L^{n-\ell},h\oplus k)$, where $(N,h)$ and $(L,k)$ are closed Einstein manifolds satisfying $\text{Ric}_h=-(\ell-1)\lambda h$ and $\text{Ric}_k=(n-\ell-1)\lambda k$, where $\lambda>0$ is constant. 
In what follows, we denote the corresponding Laplace-Beltrami operators of $(N,h)$ and $(L,k)$ by $\Delta_h$ and $\Delta_k$, respectively.
By Theorem 1.1 in \cite{CaseMalchiodi2024}, the operator $\mathscr{P}_{h\oplus k}^m$, $m<n/2$, can be factored as a product of the form
\begin{equation}\label{P_g^mFactorization}
\mathscr{P}_{h\oplus k}^m=
\begin{cases}
\prod_{s=0}^{\frac{m-2}{2}} D_{m-1-2s}, & m \quad\text{even};\\
\left(-\Delta_h -\Delta_k+\frac{(n-2)(n-2\ell)}{4}\lambda\right)\prod_{s=0}^{\frac{m-3}{2}} D_{m-1-2s}, & m \quad\text{odd,}
\end{cases}
\end{equation}
where 
\begin{equation*}\label{FactorD_t}
D_t:= \left(\Delta_h + \Delta_k\right)^2 -\frac{(n-2)(n-2\ell)}{2}\lambda \left( \Delta_h + \Delta_k \right) - 2t^2\lambda\left( \Delta_h - \Delta_h \right) + \left( \left(\frac{n-2}{2}\right)^2 - t^2 \right)\left(   \left(\frac{n-2\ell}{2}\right)^2 - t^2  \right)\lambda^2
\end{equation*}

As $(N,h)$ and $(L,k)$ are closed manifolds, the set of solutions of the corresponding closed eigenvalue problems
\[
-\Delta_h u = \alpha u\quad\text{on }\ N\quad\text{ and } \quad-\Delta_h v = \beta v\quad\text{on }\ L,
\]
are given by sequences of eigenvalues and eigenfunctions $\{(\alpha_i,u_i)\}_{i\in\mathbb{N}}$ and $\{(\beta_j,v_j)\}_{j\in\mathbb{N}}$  such that the sequences $(\alpha_i)$ and $(\beta_j)$ are nondecreasing, $ \alpha_i,\beta_j\geq0$, $u_i\in C^\infty(N)$, $v_j\in C^\infty(L)$,
\begin{equation*}\label{Eigenvalues}
-\Delta_h u_i =\alpha_i u_i,\quad -\Delta_k v_j=\beta v_j,
\end{equation*}
and $\{u_i\}_{i\in\mathbb{N}\cup\{0\}}$ and $\{v_j\}_{j\in\mathbb{N}\cup\{0\}}$ form orthonormal bases for $L_h(N)^2$ and $L_k^2(L)$ respectively (see, for instance, \cite{Chavel1984}). Define the functions
\[
e_{ij}:N\times L\to\mathbb{R}, \qquad e_{ij}(x,y):=u_i(x) v_j(y).
\]
Since the Laplace-Beltrami operator in $N\times L$ is just 
\[
\Delta_{h\oplus k} = \Delta_h + \Delta_k
\]
(see, for instance, \cite[Chapter 4]{Canzani2013}), it follows that 
\[
-\Delta_{h\oplus k}e_{ij} = (\alpha_i + \beta_j) e_{ij},
\]
and by Fubini's Theorem
\[
\int_{N\times L}e_{ij}e_{rs} dV_{h\oplus k}=\left(\int_N u_iu_r \;dV_h\right)\left(\int_Lv_jv_s \;dV_k \right) =
\begin{cases}
1 & (i,j)=(r,s)\\
0 & \text{otherwise.}
\end{cases}
\]
It follows that the set $\{e_{ij}\}_{i,j\in\mathbb{N}\cup\{0\}}$ is orthonormal in $L^2_{h\oplus k}(N\times L)$ and consists in eigenfunctions of $-\Delta_{h\oplus k}$. In fact, it is an orthonormal basis of $L^2_{h\oplus k}(N\times L)$ (see, for instance, \cite[Section II.4]{ReedSimon1980}). We next see that the operator $D_t$ is positive for suitable $t\in\mathbb{N}$. 

\begin{lemma}\label{Lemma:EigenvaluesD_t}
For $\ell\geq1$, if $n\geq 2m+2\ell-1$ and $0\leq t\leq m-1$ is an integer, then, for each $i,j\in\mathbb{N}\cup\{0\}$
\[
D_t e_{ij} = q_t(\alpha_i,\beta_j) e_{ij}
\]
where $q_t(\alpha_i,\beta_j)>0$ is a constant, depending only on $t,\alpha_i,\beta_j,\ell$ and $n$. Moreover, there constants $q(\alpha_i,\beta_j)$ are uniformly bounded from below by a constant $d_t>0$.
\end{lemma}

\begin{proof}
We prove first that $e_{ij}$ is an eigenfunction of $D_t$. For simplicity, write $A:=-\Delta_h$, $B:=-\Delta_k$, $a:=\frac{n-2}{2}$ and $b:=\frac{n-2\ell}{2}$. Then, the operator $D_t$ can be written as
\[
D_t = (A+B)^2 + 2ab\lambda(A+B) + 2t^2(A-B) + (a^2-t^2)(b^2-t^2)\lambda^2.
\]
So, given $e_{ij}=u_i v_j$, as $Au_i=\alpha_i u_i$ and $B v_j =\beta_jv_j$, the above expression yields that
\begin{equation}\label{EigenvalueD_t}
D_te_{ij}=q_t(\alpha_i,\beta_j) e_{ij}, \qquad q_t(\alpha_i,\beta_j):= (\alpha_i+\beta_j)^2 + 2\lambda(ab+t^2)\alpha_i+ 2\lambda(ab-t^2)\beta_j + (a^2-t^2)(b^2-t^2)\lambda^2.
\end{equation}
Now, notice that the hypothesis on the dimension $n$ implies that $0<b<a$. Also notice that the hypothesis on $t$ gives that $b>t\geq 0$, because $b>t$ is equivalent to the inequality $n>2t+2\ell$, and this happens because
\[
0\leq t\leq m-1 \quad \Longrightarrow\quad 2t+2\ell\leq 2m+2\ell - 2 < 2m+2\ell-1\leq n
\]
which is precisely the hypothesis on the dimension $n$. Hence, since $\ell\geq 1$, we have that $a>b>t\geq 0$ and, therefore, $a^2\geq ab\geq b^2>t^2$, so that the terms $(ab+t^2),(ab-t^2),(a^2-t^2)$ and $(b^2-t^2)$ in the expression of $q_t(\alpha_i,\beta_j)$ are all positive. As $\lambda$ is also positive, we conclude that $q(\alpha_i,\beta_j)>0$. Moreover, as $\alpha_i,\beta_j\geq 0$ for each $i,j\in\mathbb{N}\cup\{0\}$, we have that 
\begin{equation}\label{EigenvaluesUniformlyBoundedBelow}
q_t(\alpha_i,\beta_j)\geq (a^2-t^2)(b^2-t^2)\lambda^2=: d_t>0,\qquad\forall\; i,j\in\mathbb{N}\cup\{0\}
\end{equation}
and the eigenvalues $q_t(\alpha_i,\beta_j)$ are uniformly bounded from below by a positive constant. 
\end{proof}

As a previous step, we verify first that the operator $D_t$ is coercive.

\begin{lemma}\label{Lemma:OperatorD_tCoercive}
For $\ell\geq 1$, if $n\geq 2m+2\ell-1$ and $0\leq t\leq m-1$ is an integer, then the operator $D_t$ is coercive, i.e.,
\begin{equation}\label{CoercivityD_t}
\int_{N\times L}w D_t w \geq d_t \Vert w\Vert_{L_{h\oplus k}^2(N\times L)}^2,\qquad \forall\;w\in C^\infty(N\times L) . 
\end{equation}
\end{lemma}

\begin{proof}
To make the exposition clearer in this proof, relabel the eigenvalues $\{\alpha_i+\beta_j\}_{i,j}$ of the closed eigenvalue problem for $-\Delta_{h\oplus k}$ as the eigenvalues with multiplicities $\{\Lambda_i\}_{i\in\mathbb{N}}$, i.e., $\Lambda_i<\Lambda_{i+1}$ and denote the corresponding eigenspace by $H_i$, which is finite-dimensional with dimension $\sigma_i:=\dim H_i>0$. Accordingly, relabel the corresponding orthonormal basis $\{e_{ij}\}_{i,j}$  for $L^2_{h\oplus k}(N\times L)$ as $\{E_i^j\;:\; i\in\mathbb{N}, j\in\{1,\ldots \sigma_i\}\}$, so that, in this new notation,
\[
-\Delta_{h\oplus k}E_{i}^j = \Lambda_i E_i^j,\qquad\text{ for every }i\in\mathbb{N}\quad\text{and every }\quad j=1,\ldots\sigma_i,
\]
and $H_i=\text{span}\{E_i^1,\ldots E_i^{\sigma_i}\}$. Hence,
\[
L_{h\oplus k}^2(N\times L) = \bigoplus_{i=1}^\infty H_i
\]
in the sense that any $w\in L_{h\oplus k}^2(N\times L)$ can be written as
\[
w = \sum_{i=1}^\infty w_i,\qquad w_i\in H_i.
\]
Further, relabel the eigenvalues $\{q(\alpha_i,\beta_j)\}_{i,j\in\mathbb{N}\cup\{0\}}$ of $D_t$ as $\{Q_i^j\;:\; i\in \mathbb{N}, j\in\{1,\ldots\sigma_i\}\}$, so that \eqref{EigenvalueD_t} and \eqref{EigenvaluesUniformlyBoundedBelow} can be rewritten as
\begin{equation*}\label{RelabeledEigenvaluesD_t}
D_t E_{i}^j = Q_i^jE_{i}^j,\qquad Q_{i}^j\geq d_t,\qquad \text{ for each }i\in\mathbb{N} \text{ and each }j=1,\ldots\sigma_i.
\end{equation*}
We next verify \eqref{CoercivityD_t} in several steps. 

\smallskip
\textbf{Step 1.} \emph{The inequality \eqref{CoercivityD_t} is true for every $w\in H_i$ and $D_tw\in H_i$.}

\smallskip

First, if $w\in H_i$  for some $i\in\mathbb{N}$, then
\[
w = \sum_{j=1}^{\sigma_i}c_{ij}E_i^j,\qquad \Vert w\Vert_{L_{h\oplus k}^2(N\times L)} =  \sum_{j=1}^{\sigma_i}\vert c_{ij}\vert^2 
\]
where $c_{ij}\in\mathbb{R}$. Since the eigenfunctions of the closed eigenvalue problem are smooth (see \cite{Chavel1984}) then $E_i^j$ are smooth and $w\in C^\infty(N\times L)$. Hence,
\begin{equation}\label{H_iInvariantUnderD_t}
D_t w = \sum_{j=1}^{\sigma_j} D_t c_{ij}E_i^j = \sum_{j=1}c_{ij}Q_i^jE_i^j\in H_i,
\end{equation}
and therefore,
\begin{align*}
\int_{N\times L}w D_t  w \;dV_{h\oplus k} 
&= \sum_{j=1}^{\sigma_i} \sum_{J=1}^{\sigma_i}\int_{N\times L} c_{ij}c_{iJ} E_i^j D_t E_i^J \;dV_{h\oplus k} \\
&=  \sum_{j=1}^{\sigma_i} \sum_{J=1}^{\sigma_i}c_{ij}c_{iJ} Q_j^J\int_{N\times L} E_i^j  E_i^J \;dV_{h\oplus k} \\
&= \sum_{j=1}^{\sigma_i}\vert  c_{ij}\vert^2Q_i^j 
\geq \sum_{j=1}^{\sigma_i}\vert  c_{ij}\vert^2d_t = d_t\Vert w\Vert_{L_{h\oplus k}^2(N\times L)}^2
\end{align*}
which is \eqref{CoercivityD_t} in this case.

\smallskip
\textbf{Step 2.} \emph{For any $J\in\mathbb{N}$, the inequality \eqref{CoercivityD_t} is true for every $w\in \bigoplus_{i=1}^JH_i$.}

\smallskip

Take $J\in\mathbb{N}$ and  $w\in \bigoplus_{i=1}^JH_i$, then 
\[
w=\sum_{i=1}^J w_i, \quad w_i\in H_i,\quad\text{and}\quad \Vert w\Vert_{L_{h\oplus k}^2(N\times L)}^2 = \sum_{i=1}^J \Vert w_i\Vert_{L_{h\oplus k}^2(N\times L)},
\]
because the eigenspaces are mutually orthogonal in $L_{h\oplus k}^2(N\times L)$. Moreover, as each $w_i$ is smooth, then also $w\in C^\infty(N\times L)$ is smooth. Since $D_tw_i\in H_i$ by \eqref{H_iInvariantUnderD_t} and as the eigenspaces are mutually orthogonal, since
the  inequality \eqref{CoercivityD_t} is valid for $w_i\in H_i$, then
\[
\int_{N\times L} w D_t w \;dV_{h\oplus k} = \sum_{i=1}^J \sum_{j=1}^J \int_{N\times L} w_iD_t w_j\;dV_{h\oplus k} =\sum_{i=1}^J\int_{N\times L}w_iD_t w_i \;dV_{h\oplus k}\geq \sum_{i=1}^J d_t\Vert w_i\Vert_{L_{h\oplus k}^2(N\times L)}^2 = d_t\Vert w\Vert_{L_{h\oplus k}^2(N\times L)}^2,
\]
as we wanted.

\smallskip
\textbf{Step 3.} \emph{For any $w\in C^\infty(N\times L)$, the inequality \eqref{CoercivityD_t} is true.}

\smallskip

Take $w\in C^\infty(N\times L)$ and consider $w_i\in H_i$ for each $i\in\mathbb{N}$ such that $W_J:=\sum_{i=1}^J w_i\to w$ in $L_{h\oplus k}^2(N\times L)$. We next see that $D_tW_J\to D_t w$ in $L_{h\oplus k}^2(N\times L)$.   Since $D_t$ is a fourth order linear elliptic operator with domain $C^{\infty}(N\times L)$, it can be extended to a bounded linear operator $D_t:H^4_{h\oplus k}(N\times L)\to L^2_{h\oplus k}(N\times L)$ (see, for instance, \cite[Proposition 2.1]{LawsonMichelsohn1989}).  Now, write each $w_i\in H_i$ in terms of the orthonormal basis $E_i^j$ as
\[
w_i = \sum_{j=1}^{\sigma_i}c_{ij}E_i^j,
\]
so that
\[
w = \sum_{i=1}^\infty\sum_{j=1}^{\sigma_i}c_{ij}E_i^j, \qquad \Vert w\Vert_{L_{h\oplus k}^2(N\times L)}^2=\sum_{i=1}^\infty\Vert w_i\Vert_{L_{h\oplus k}^2(N\times L)}^2 = \sum_{i=1}^\infty\sum_{j=1}^{\sigma_j}\vert c_{ij}\vert^2,
\]
and recall that for any $K\in\mathbb{N}$ the norm in $H_{h\oplus k}^K(N\times L)$ is equivalent to the spectral norm
\begin{equation}\label{SpectralNorms}
\Vert w\Vert^2_{spec_K}:= \sum_{i=1}^\infty\sum_{j=1}^{\sigma_i}(1+\Lambda_i^K)\vert c_{ij}\vert^2<\infty
\end{equation}
(see \cite[Proposition 56]{Canzani2013} or \cite[Section 7.3]{LionsMagenes1972}). Since $C^\infty(N\times L)\subset H_{h\oplus k}^4(N\times L)$, then $w, W_J,w-W_J\in H_{h\oplus k}^4(N\times L)$ and 
\[
\Vert w-W_J\Vert_{spec_4} = \sum_{i=J+1}^\infty\sum_{j=1}^{\sigma_i} (1+\Lambda_i^4)\vert c_{ij}\vert^2 \to 0
\]
because $w\in H_{h\oplus k}^4(N\times L)$ and the above sum is the tail of the series given in the spectral norm of $w$ in \eqref{SpectralNorms} for $K=4$. Then, $W_J\to w$ in $H_{h\oplus k}^4(N\times L)$ and the continuity of $D_t$ yields that
\[
D_tW_J=\sum_{i=1}^JD_tw_i \to D_tw\qquad \text{in } \; L_{h\oplus k}^2(N\times L).
\]
Since also $W_J\to w$ in $L_{h\oplus k}^2(N\times L)$, then the sequence $(W_J)_J$ is bounded in $L_{h\oplus k}^2(N\times L)$ and $W_J\rightharpoonup w$ weakly in $L_{h\oplus k}^2(N\times L)$. Hence, as $w_i,D_t w_i\in H_i$ by \eqref{H_iInvariantUnderD_t}, and as the spaces $H_i$ are mutually orthogonal, by Step 1 we have that
\begin{align*}
\int_{N\times L} wD_t w \;dV_{h\oplus k} 
&= \lim_{I,J\to\infty} \sum_{i=1}^I\sum_{j=1}^J \int_{N\times L}w_i D_t w_j \;dV_{h\oplus k} \\
&=\lim_{J\to\infty}\sum_{i=1}^J\int_{L\times N} w_i D_t w_i \;dV_{h\oplus k}\\
&\geq \lim_{J\to\infty}\sum_{i=1}^J d_t \Vert w_i\Vert_{L_{h\oplus k}^2(N\times L)}^2 = d_t \Vert w\Vert_{L_{h\oplus k}^2(N\times L)}^2,
\end{align*}
and we end the proof of the result.
\end{proof}

We can now prove that, under suitable dimensions, the $GJMS$-operator on special Einstein products is coercive. Denote by
\[
\lambda_1^m(N\times L,h\oplus k):= \inf_{w\in C^\infty(N\times L)\smallsetminus\{0\}} \frac{\int_{N\times L}w\mathscr{P}_{h\oplus  k}^mw}{\Vert w\Vert_{L_{h\oplus k}^2(N\times L)}^2}
\]
the first eigenvalue of $\mathscr{P}_{h\oplus k}^m$.
\begin{proposition}\label{Proposition:P_hkCoercive}
If $n\geq 2m+2\ell -1$, then $\lambda_1^m(N\times L,h\oplus k)>0$. In particular, $\mathscr{P}_{h\oplus k}^m$ is coercive.
\end{proposition}

\begin{proof}
First, notice that for any $0\leq s,t\leq m-1$ and any $i,j\in\mathbb{N}\cup\{0\}$, we have that
\begin{equation}\label{CommutativityD_tD_s}
D_tD_s e_{ij} =D_t q_s(\alpha_i,\beta_j)e_{ij}= q_t(\alpha_i,\beta_j) q_s(\alpha_i,\beta_j) e_{ij},
\end{equation}
where $q_s(\alpha_i,\beta_j)\geq d_s$ and $q_t(\alpha_i,\beta_j)\geq d_t$ for any $i,j\in\mathbb{N}\cup\{0\}$ by Lemma \ref{Lemma:EigenvaluesD_t}. Recall now the decomposition of the operator $\mathscr{P}_{h\oplus k}^m$ given in \eqref{P_g^mFactorization}, which depends on the parity of $m$. We next prove  that there exists a constant $d_0>0$ such that
\begin{equation}\label{CoercivityOperatorEinsteinSpecialProducts}
\int_{N\times L}w \mathscr{P}_{h\oplus k}^m w\;dV_{h\oplus k} \geq d_0 \Vert w\vert_{L_{h\oplus k}^2(N\times L)}^2,\qquad\forall\; w\in C^\infty(N\times L),
\end{equation}
which immediately implies that $\lambda_{1}^m(N\times,h\oplus k)\geq d_0$ and the coercivity of $\mathscr{P}_{h\oplus k}^m$. We prove this for $m$ odd, the case of $m$ even being similar and simpler.

The proof of this fact proceed along the same lines as in the proof of Lemma \eqref{Lemma:OperatorD_tCoercive} and, consequently, it suffices to show inequality \eqref{CoercivityOperatorEinsteinSpecialProducts} for any finite linear combination of the Hilbert basis $\{e_{ij}\}$.

First, let $t_s:=m-1-2s$ for any integer $0\leq \frac{m-3}{2}$. Then, given any $i,j\in\mathbb{N}\cup\{0\}$, by \eqref{CommutativityD_tD_s}, we have that
\begin{align*}
\mathscr{P}_{h\oplus k}^m e_{ij} &= \left(-\Delta_h-\Delta_k + \frac{(n-2)(n-2\ell)}{4}\lambda\right) \prod_{s=0}^{\frac{m-3}{2}}D_{t_s}e_{ij}\\
&= \prod_{s=0}^{\frac{m-3}{2}}q_{t_s}(\alpha_i,\beta_j)\left(-\Delta_h-\Delta_k + \frac{(n-2)(n-2\ell)}{4}\lambda\right) e_{ij}\\
&= \left[\left(\alpha_i+\beta_j+\frac{(n-2)(n-2\ell)}{4}\lambda\right)\prod_{s=0}^{\frac{m-3}{2}}q_{t_s}(\alpha_i,\beta_j)\right]e_{ij},
\end{align*}
i.e., the Hilbert basis consists of eigenvectors of $\mathscr{P}_{h\oplus k}^m$ with eigenvalues 
\[
q_{ij}:=\left(\alpha_i+\beta_j+\frac{(n-2)(n-2\ell)}{4}\lambda\right)\prod_{s=0}^{\frac{m-3}{2}}q_{t_s}(\alpha_i,\beta_j).
\]

Since $n\geq 2m+2\ell-1$, then $2< t_s\leq m-1$ for every integer $0\leq s\leq \frac{m-3}{2}$, and Lemma \ref{Lemma:EigenvaluesD_t} gives a lower bound $d_{t_s}>0$ for $q_{t_s}(\alpha_i,\beta_i)$. Hence,
\begin{equation}\label{BoundEigenvaluesP_hk}
q_{ij}\geq \left(\alpha_i+\beta_j+\frac{(n-2)(n-2\ell)}{4}\lambda\right)\prod_{s=0}^{\frac{m-3}{2}}d_{t_s}\geq \frac{(n-2)(n-2\ell)}{4}\lambda\prod_{s=0}^{\frac{m-3}{2}}q_{t_s}=:d_0>0,\qquad\forall\; i,j\in\mathbb{N}\cup\{0\},
\end{equation}
because the eigenvalues $\alpha_i,\beta_j\geq 0$. Hence, given $w\in\text{span}\{e_{ij}\;:\; i,j\in\mathbb{N}\cup\{0\} \}$, $w=\sum_{i=0}^I\sum_{j=0}^J c_{ij}e_{ij}$, for some $I,J\in\mathbb{N}\cup\{0\}$ we have that
\[
w\in C^{\infty}(N\times L),\qquad \Vert w\Vert_{L_{h\oplus k}^2(N\times L)}^2 = \sum_{i=0}^I\sum_{j=0}^J \vert c_{ij}\vert^2,
\]
and by \eqref{BoundEigenvaluesP_hk} we deduce that
\begin{align*}
\int_{N\times L} w\mathscr{P}_{h\oplus k}^m w \;dV_{h\oplus k} 
&= \sum_{i=0}^{I}\sum_{j=0}^J\sum_{i'=0}^{I}\sum_{j'=0}^{J}c_{ij}c_{i'j'} \int_{N\times L} e_{i'j'}\mathscr{P}_{h\oplus k}^m e_{ij} \;dV_{h\oplus k}
 =\sum_{i=0}^{I}\sum_{j=0}^J\sum_{i'=0}^{I}\sum_{j'=0}^{J}c_{ij}c_{i'j'} \int_{N\times L} e_{i'j'}q_{ij} e_{ij}\;dV_{h\oplus k}\\
&= \sum_{i=0}^{I}\sum_{j=0}^J\sum_{i'=0}^{I}\sum_{j'=0}^{J}c_{ij}c_{i'j'} q_{ij}\int_{N\times L} e_{i'j'} e_{ij}\;dV_{h\oplus k}
=\sum_{i=0}^{I}\sum_{j=0}^J\vert c_{ij}\vert^2 q_{ij} \\
&\geq d_0 \sum_{i=0}^{I}\sum_{j=0}^J\vert c_{ij}\vert^2 = d_0 \Vert w\Vert_{L_{h\oplus k}^2(N\times L)}^2.
\end{align*}
The rest of the proof follows as in Step 2 and Step 3 in Lemma \ref{Lemma:OperatorD_tCoercive}, using in this case that $\mathscr{P}_{h\oplus k}^m$ is a $2m$-order linear elliptic operator with domain $C^\infty(N\times L)$ that can be extended to a bounded linear operator $\mathscr{P}_{h\oplus k}^m: H^{2m}_{h\oplus k}(N\times L)\to L_{h\oplus k}^2(N\times L)$.
\end{proof}

We can now prove Corollary \ref{Corllary:SpecialEinsteinProducts}.

\begin{proof}[Proof of Corollary \ref{Corllary:SpecialEinsteinProducts}]
Let $n\geq D(m,\ell)$. By Theorem 2 in \cite{CaseMalchiodi2024}, the operator $\mathscr{P}_{h\oplus k}^m$ satisfies property \eqref{Property:PPP}, and by \eqref{ConstantD(m,l)}, $n\geq 2m+2\ell-1$, so Proposition \ref{Proposition:P_hkCoercive} yields that the operator is also coercive. Finally, to apply Theorem \ref{Theorem:Main}, we need to check condition \eqref{OrbitsPositiveDimensional}: Indeed, $\Gamma:=\Gamma_N\times\Gamma_L$ is a closed subgroup of $\text{Isom}(N\times L,h\oplus k)$ and the $\Gamma$-orbits are given by
\[
\Gamma(x,y) = \Gamma_Nx\times \Gamma_L y,\qquad\forall\;(x,y)\in N\times L.
\]
Therefore 
\[
\dim\Gamma(x,y) =\dim\Gamma_Nx + \dim \Gamma_Ly
\]
and \eqref{OrbitsPositiveDimensional} follows directly from \eqref{OrbitsSpecialEinsteinProduct}. Hence, all the hypotheses in Theorem \ref{Theorem:Main} are satisfied and problem \eqref{MainEquation:Q-curvatureEquation} admits an infinite number of $\Gamma$-invariant solutions, one of them positive and the rest change sign.
\end{proof}

\bibliographystyle{siam}
\bibliography{CombinedReferences}

\begin{flushleft}

\textbf{Juan Carlos Fernández}\\

Departamento de Matemáticas, Facultad de Ciencias\\
Universidad Nacional Autónoma de México\\
Circuito Exterior, Ciudad Universitaria\\
04510 Coyoacán, Ciudad de México, Mexico\\
\texttt{jcfmor@ciencias.unam.mx}

\end{flushleft}

\end{document}